\documentclass[10pt, twocolumns]{IEEEtran}

\usepackage{cite}
\usepackage{amsmath,amssymb,amsfonts}
\usepackage{algorithmic}
\usepackage{textcomp}
\usepackage{color}
\usepackage{bm}
\usepackage{graphicx}
\usepackage{enumerate}
\usepackage{multicol}

\newtheorem{lemma}{Lemma}
\newtheorem{assumption}{Assumption}

\newtheorem{remark}{Remark}
\newtheorem{corollary}{Corollary}
\newtheorem{theorem}{Theorem}
\newtheorem{example}{Example}

\title{
\huge{
Primal-dual Accelerated Mirror-Descent Method 
for Constrained Bilinear Saddle-Point Problems
}
}

\author{Weijian~Li, Xianlin~Zeng, and Lacra Pavel

\thanks{W.~Li and L. Pavel are with the Department of Electrical
and Computer Engineering, University of Toronto, Toronto, ON,
M5S 3G4, Canada. E-mails: 
\texttt{weijian.li@utoronto.ca},
\texttt{pavel@control.utoronto.ca}.
}

\thanks{X.~Zeng is with the Key Laboratory of Intelligent Control and
	Decision of Complex Systems, School of Automation, Beijing Institute
	of Technology, Beijing 100081, China. E-mail: \texttt{xianlin.zeng@bit.edu.cn}.}

}

\begin{document}

\maketitle

\begin{abstract}
We develop a first-order accelerated algorithm for a class of constrained bilinear saddle-point problems with applications to network systems.
The algorithm is a modified time-varying primal-dual version of an accelerated mirror-descent dynamics.
It deals with constraints such as simplices and convex set constraints effectively, and converges with a rate of $\mathcal O(1/t^2)$.
Furthermore, we employ the acceleration scheme to constrained distributed optimization and bilinear zero-sum games, and obtain two variants of distributed accelerated algorithms.
\end{abstract}


\section{INTRODUCTION}

Constrained bilinear saddle-point problems (BSPPs) aim at seeking saddle points that minimize {\color{blue}global functions} in one direction and maximize them in the other, where decision variables are limited to certain feasible regions or ranges.
BSPPs arise from convex optimization with affine constraints via Lagrangian methods \cite{ruszczynski2011nonlinear}, 
while constraints are significant in applications  including multi-robot motion planning, resource allocation in communication networks, and economic dispatch in power grids. 
Constrained BSPPs appear in various scenarios such as  robust optimization, signal processing, and machine learning \cite{kovalev2022accelerated}, and thus,
a substantial  effort has been paid to design effective algorithms to solve them.
Particularly, first-order methods are frequently used since they only involve gradient information and are suitable for large-scale problems.
{\color{blue}First-order primal-dual methods have been widely investigated}. They are powerful for constrained optimization/game problems, and easy to be implemented both in centralized and distributed manners \cite{arrow1958studies, cherukuri2017saddle, cherukuri2017role}.

Motivated by big-data advances and privacy concerns encountered in practice, the areas of distributed optimization and games over multi-agent systems have received much attention, with applications to unmanned systems, smart grids, and social networks \cite{nedic2010constrained, gharesifard2013distributed}. Some of the formulations could be cast into constrained BSPPs, and then, solved by primal-dual methods \cite{gharesifard2013distributed, zeng2018distributed, lei2016primal, liang2019distributed}.
For instance, a discrete-time primal-dual algorithm was designed for a constrained distributed optimization problem in \cite{lei2016primal}, while a continuous-time version was proposed for a distributed extended monotropic optimization problem in \cite{zeng2018distributed}.
In \cite{liang2019distributed}, a distributed primal-dual algorithm with a fixed stepsize was developed to handle coupled constraints. Additionally, a primal-dual framework was explored for two-network zero-sum games in \cite{gharesifard2013distributed}.
However, existing results indicate that standard primal-dual methods could only achieve asymptotic convergence with a rate of $\mathcal O(1/t)$ for general convex cost functions \cite{zeng2018distributed, lei2016primal}.

Various accelerated protocols have been explored for convex optimization, for example the well-known heavy-ball algorithm \cite{polyak1987introduction}.
The author of \cite{nesterov2003introductory}  showed that first-order methods have an optimal rate of $\mathcal O(1/t^2)$ in the worst case.
Then the Nesterov's accelerated scheme was proposed in \cite{nesterov2003introductory}, and a continuous-time 
version was developed in \cite{su2014differential}.
{\color{blue}
In \cite{scieur2020regularized}, an accelerated method was explored under regularity conditions.
In \cite{poveda2021robust}, a class of zero-order optimization dynamics was proposed with acceleration and restarting mechanisms.}
However, these results are mainly devoted to unconstrained optimization.
For consensus-based distributed optimization,
a fast gradient algorithm, with a rate of $\mathcal O(1/t^{1.4-\epsilon})$, was designed in \cite{jakovetic2014fast},
while a Nesterov's accelerated algorithm, with a rate of $\mathcal O(1/t^{2 - \epsilon})$, was developed in \cite{qu2019accelerated}.
However, both of them were primal-based and only handled consensus type constraints.
In order to solve convex optimization with affine constraints or BSPPs,
attempts have been made by combining accelerated protocols with  primal-dual methods. 
For instance, an accelerated primal-dual method, using adaptive parameters, was designed with a rate of $\mathcal O(1/t^2)$ in \cite{xu2017accelerated}.
In \cite{kovalev2022accelerated}, an accelerated gradient method, with an optimal linear rate, was developed.
In \cite{salim2022optimal}, lower bounds of the gradient computation burden were provided, and then, a corresponding  algorithm was proposed with complexity guarantees.
However, the algorithms in \cite{xu2017accelerated, kovalev2022accelerated, salim2022optimal} require strong convexity on cost functions, which does not hold in applications such as data analysis and machine learning \cite{necoara2019linear}.
In  \cite{zeng2022dynamical}, by incorporating the Nesterov's acceleration into a primal-dual framework, a continuous-time accelerated dynamics was designed for strictly convex cost functions.
Afterwards, the protocol was adopted to seek Nash equlibria for two-network bilinear zero-sum games in \cite{zeng2023distributed}.
It has been shown that the algorithms in \cite{zeng2022dynamical, zeng2023distributed} achieved a rate of $\mathcal O(1/t^2)$, but neither of them can handle optimization/game problems beyond affine constraints, including set constraints.
To the best of our knowledge, few accelerated primal-dual algorithms, dealing with strictly convex cost functions and set constraints at the same time, achieve a rate of $\mathcal O(1/t^2)$.

The mirror-descent algorithm,  proposed in \cite{nemirovskij1983problem}, is a generalization of gradient-based methods by introducing a Bregman distance function in place of the Euclidean distance.
It can solve constrained optimization problems  \cite{aharon2001order}, and has generated a lot of interest. For instance, linear rates of mirror-descent and actor-critic dynamics were analyzed in \cite{gao2020continuous, gao2022continuous}.
In \cite{gao2023second}, a second-order mirror-descent dynamics was proposed for merely monotone games.
In \cite{chen2021distributed}, a primal-dual version was designed for distributed optimization.
In \cite{krichene2015accelerated}, an accelerated mirror-descent scheme, regarded as an extension of the Nesterov's acceleration, was developed with a rate of $\mathcal O(1/t^2)$.
However, existing literature has not provided any accelerated mirror-descent algorithm under a primal-dual framework.

In this paper, we focus on developing such a continuous-time algorithm.
{\color{blue}On the one hand, continuous-time algorithms have received a flurry of research interest with applications in multi-agent systems.
For cyber-physical systems with physical systems acting or moving in real time, they are quite natural, and meanwhile, can be implemented by hardware devices such as analog circuits \cite{forti2004generalized}. The continuous-time design has been reported in existing literature, including \cite{arrow1958studies, cherukuri2017saddle, liang2019distributed, gharesifard2013distributed, yang2016multi, poveda2021robust}. On the other hand, continuous-time algorithms can serve as a tool for understanding, analyzing and generalizing the discrete-time algorithms \cite{luo2024universal, luo2022differential, su2014differential}. }
Our main contributions are summarized as follows.

\begin{enumerate}[a)]
\item We propose a time-varying primal-dual mirror-descent algorithm to solve a constrained BSPP, which covers many problems including convex optimization with affine constraints \cite{zeng2022dynamical, xu2017accelerated}, constrained distributed optimization \cite{lei2016primal}, and bilinear games \cite{zeng2023distributed}.
The algorithm is a modified primal-dual version of the dynamics in
\cite{krichene2015accelerated}, has a unique solution, and can efficiently handle set constraints, simplices, etc.

\item We establish that the algorithm converges with a rate of $\mathcal O(1/t^2)$ (in terms of the duality gap function), which is faster than the primal-dual methods in \cite{arrow1958studies, cherukuri2017saddle, cherukuri2017role}. 

\item Resorting to the accelerated protocol, we provide two variants of distributed algorithms for constrained distributed optimization and for bilinear zero-sum games.
We show that they have faster rates than those of \cite{zeng2018distributed, lei2016primal, yang2016multi, gharesifard2013distributed}.
\end{enumerate}

This paper is organized as follows. Some preliminaries are introduced in Section \uppercase\expandafter{\romannumeral 2}. 
Section \uppercase\expandafter{\romannumeral 3} formulates the problem, and presents the accelerated algorithm.
Our main results are provided in Section \uppercase\expandafter{\romannumeral 4}, and two distributed algorithms  are obtained for network optimization and games in Section \uppercase\expandafter{\romannumeral 5}. Finally, concluding remarks are given in Section
\uppercase\expandafter{\romannumeral 6}.

\section{Preliminary Background}

Let $\mathbb R$,  $\mathbb N$, $\mathbb R^p$ and $\mathbb R^{p \times q}$ be the set of real numbers, nonnegative integers, $p$-dimensional real column vectors and $p$-by-$q$ real matrices. 
Denote $0_p$ as the $p$-dimensional zero vector, and $I_p$ as the $p$-by-$p$ identity matrix.
Denote $(\cdot)^{\rm T}$, $\otimes$ and $\Vert \!\cdot\! \Vert$ as the transpose, the Kronecker product, and the Euclidean norm.
For $x, y \in \mathbb R^p$, their Euclidean inner product is  $x^{\rm T} y$ or $\langle x, y \rangle$.
Let $X \times Y$ be the Cartesian product of sets $X$ and $Y$.
For a differentiable function $S(x, y)$, $\nabla_x S(x, y)$ is the partial gradient of $S$ with respect to $x$.
Let $\{x_n(t)\}_{n \in \mathbb N}$ be a family of functions from $\mathbb R$ to $\mathbb R^p$.
The sequence is equi-continuous 
if for any $\epsilon > 0$, there always exists $\delta > 0$ such that
$\Vert x_n(t + h) - x_n(t) \Vert < \epsilon,  \forall |h| < \delta, \forall n \in \mathbb N$.

Let $\Omega \subset \mathbb R^p$ be a {\color{blue} closed convex set}.
The normal cone of $\Omega$ at $x \in \Omega$ is 
$\mathcal N_{\Omega}(x) =
\{v \in \mathbb R^p | v^{\rm T}(y-x) \le 0, \forall y \in \Omega\}.$
Take
$P_{\Omega}(u) = {\rm argmin}_{v \in \Omega} \Vert v - u\Vert$.
Then
$\langle v - P_{\Omega}(u), u - P_{\Omega}(u)\rangle \le 0,
~\forall u \in \mathbb R^p, ~\forall v \in \Omega$.

{\color{blue}
Take $\mathbb E = \mathbb R^p$.
Given an extended real-valued function $\psi: \mathbb E \rightarrow [- \infty, + \infty]$, it is proper if it does not attain the value $- \infty$, and there exists at least one $x \in \mathbb E$ such that $\psi(x) < \infty$.
Besides, it is closed if its epigraph is closed.
The domain of $\psi$ is ${\rm dom}(\psi) = \{x \in \mathbb E: \psi(x) < \infty\}$.
The function $\psi^* :  \mathbb E^* \rightarrow [-\infty, +\infty]$, defined by $\psi^*(z) = \sup_{x\in \mathbb E} \{x^T z - \psi(x)\}$, is called the conjugate function of $\psi$, where $\mathbb E^*$ is the dual space of $\mathbb E$. If $\psi$ is proper, then $\psi^*$ is closed and convex.
The Bregman divergence of a proper, closed and convex function $\psi$ is 
$D_{\psi}(x, y) = \psi(x) - \psi(y) - \langle x - y, 
\nabla \psi(y)\rangle$.
}
Referring to \cite{diakonikolas2019approximate}, we have the following result.
\begin{lemma}
\label{lem:gen}
Let  $\psi: \mathbb R^p \rightarrow \mathbb R \cup \{\infty\}$ be a closed, proper and strongly convex function, 
where {\color{blue}$\Omega := {\rm dom}(\psi)$ is a closed convex set}. Then $\psi^*$ is  differentiable and convex such that
\begin{equation}
\label{conjugate}
\nabla {\psi}^*(z) = {\rm argmin}_{x \in \Omega} \{-x^{\rm T} z + \psi(x)\}.
\end{equation}
\end{lemma} 

Consider a multi-agent network described by an undirected graph $\mathcal G(\mathcal V, \mathcal E)$, where $\mathcal V = \{1, \dots, n\}$ is the node set and $\mathcal E \subset \mathcal V \times \mathcal V$ is the edge set. The graph $\mathcal G$ is associated with an adjacency matrix $\mathcal A = [a_{ij}] \in \mathbb R^{n \times n}$ such that $a_{ij} = a_{ji} > 0$ if $(i, j) \in \mathcal E$, and $a_{ij} = 0$ otherwise.
The Laplacian matrix $\mathcal L$ is $\mathcal L = \mathcal D - \mathcal A$, where $\mathcal D = {\rm diag}\{d_i\}\in \mathbb R^{n \times n}$, and $d_i = \sum_{j \in \mathcal V} a_{ij}$.
Graph $\mathcal G$ is connected if there exists a path between any pair of distinct nodes.

\section{PROBLEM SETUP AND ALGORITHM}

In this section, we give the problem statement and  propose an accelerated algorithm.

\subsection{Problem Formulation}

Consider a constrained BSPP given by
\begin{equation}
\begin{aligned}
\label{form}
\min_{x \in X}~\max_{y \in Y}~
L(x,y) := F(x) + x^{\rm T} H y - G(y)
\end{aligned}
\end{equation}
where $X \subset \mathbb R^p$ and $Y \subset \mathbb R^q$  are  closed constraint sets, $H \in \mathbb R^{p \times q}$ is a known matrix, and moreover, $F: X \rightarrow \mathbb R$ and $G: Y \rightarrow \mathbb R$ are continuous functions.
A profile $(x^\star, y^\star) \in X \times Y$ is called a saddle point (or a solution) to (\ref{form}) if 
$$
L(x^\star,y) \le L(x^\star,y^\star) \le L(x,y^\star),
~~\forall x \in X, ~\forall y \in Y.
$$

{\color{blue}We make the following well-known assumption for (\ref{form})}. 
\begin{assumption}
\label{ass:convex}
{\color{blue} The sets $X$ and $Y$} are nonempty, closed and convex.
{\color{blue} The functions $F$ and $G$} are continuously differentiable and strictly convex on two open sets containing $X$ and $Y$.
{\color{blue} The gradients of $F$ and $G$} are Lipschitz continuous over $X$ and $Y$.
{\color{blue}Moreover, there exists a unique saddle point to (\ref{form}).}
\end{assumption}

\begin{remark}
Problem (\ref{form}) is a general model, and appears in constrained optimization, bilinear games and multi-agent systems \cite{ruszczynski2011nonlinear, lei2016primal, gharesifard2013distributed}. 
{\color{blue}
It generalizes the models in 
\cite{zeng2022dynamical, zeng2023distributed} by allowing the constraints $X$ and $Y$, and covers the problems in \cite{kovalev2022accelerated, xu2017accelerated, salim2022optimal} since $F$ and $G$ are only required to be strictly convex rather than strongly convex.}
\end{remark}

First-order primal-dual methods in \cite{arrow1958studies, cherukuri2017saddle, cherukuri2017role}  can be employed to solve (\ref{form}). 
However, they suffer from a slow convergence rate.
Existing works indicate that mirror-descent methods, powerful for constrained optimization, can achieve acceleration \cite{nemirovskij1983problem, aharon2001order, krichene2015accelerated}.
The observations motivate us to design a novel accelerated algorithm by incorporating an accelerated mirror-descent protocol into a primal-dual framework.

\subsection{Accelerated Algorithm Design}

Consider $X = \mathbb R^p$ and $Y = \mathbb R^q$.
Referring to \cite{arrow1958studies, cherukuri2017saddle, cherukuri2017role}, a first-order primal-dual dynamics for (\ref{form}) is given by
\begin{equation}
\begin{aligned}
\label{primdual}
\dot x =& - \nabla_x L(x, y) =
- \nabla F(x)  - Hy  \\
\dot y =&~~~ \nabla_y L(x, y) =
-\nabla G(y) + H^{\rm T} x
\end{aligned}	
\end{equation}
where $x(0) = x_0 \in \mathbb R^p$, $y(0) = y_0 \in \mathbb R^q$, and  $\dot x = {\rm d}x(t)/{\rm d}t$ is the derivative of $x$ with respect to time $t$. We omit time $t$ in remaining of this paper without causing confusions.
Clearly, the variable $x$ descents in the direction of negative gradient of $L$, while $y$ ascents in the other direction.

Given a closed convex set $\Omega \subset \mathbb R^p$ and a continuously differential and convex function $f: \Omega \rightarrow \mathbb R$, an accelerated mirror-descent algorithm is designed in \cite{krichene2015accelerated} as
\begin{equation}
\begin{aligned}
\label{acce:mirror}
\dot x = \frac rt\big[\nabla \psi^*(z) - x\big],
~\dot z = -\frac tr \nabla f(x)
\end{aligned}
\end{equation}
where $r > 0$, $\psi$ is a generating function with ${\rm dom}(\psi) = \Omega$, $\psi^*$ is the convex conjugate of $\psi$, $z(0) = z_0$ and $x(0) = x_0 = \nabla \psi^*(z_0)$.
If $\nabla \psi^*(z) = z$, then (\ref{acce:mirror}) recovers the continuous-time Nesterov's accelerated algorithm in \cite{su2014differential}.
By \cite[Th. 2]{krichene2015accelerated}, (\ref{acce:mirror}) converges with a rate of $\mathcal O(1/t^2)$ if $r \ge 2$.

Inspired by (\ref{primdual}) and (\ref{acce:mirror}), we propose a  continuous-time primal-dual accelerated algorithm as 
\begin{equation}
\begin{aligned}
\label{alg}
\dot x =& \frac rt \big[\nabla\psi^*(u) - x  \big], ~
\dot y = \frac rt \big[\nabla\phi^*(v) - y \big] \\
\dot u =& \frac tr \big[\! -\nabla F(x) - H(y+ \frac tr \dot y) + \nabla \psi(x + \frac t r \dot x)- u\big] \\
\dot v =& \frac tr \big[\! -\nabla G(y) + H^{\rm T} (x + \frac tr \dot x) + \nabla\phi(y + \frac t r \dot y) - v \big]
\end{aligned}
\end{equation}
where {\color{blue}$\psi: \mathbb R^p \rightarrow \mathbb R \cup \{\infty\}$ and $\phi: \mathbb R^q \rightarrow \mathbb R \cup \{\infty\}$
are proper, closed and convex functions with 
${\rm dom}(\psi) = X$ and ${\rm dom}(\phi) = Y$,}
$x(0) = x_0 = \nabla \psi^*(u_0)$, 
$u(0) = u_0 = \nabla \psi(x_0)$,
$y(0) = y_0 = \nabla \phi^*(v_0)$, and 
$v(0) = v_0 = \nabla \phi(y_0)$.

The following assumption is needed for $\psi$ and $\phi$
\cite{krichene2015accelerated, gao2022continuous}.
\begin{assumption}
\label{ass:genlip}
{\color{blue}
Both $\psi: \mathbb R^p \rightarrow \mathbb R \cup \{\infty\}$ and $\phi: \mathbb R^q \rightarrow \mathbb R \cup \{\infty\}$
are proper, closed, and strongly convex functions, where ${\rm dom}(\psi)$ and ${\rm dom}(\phi)$ are nonempty, closed and convex sets, which equal to $X$ and $Y$, respectively.
They are also continuously differentiable with Lipschitz continuous gradients on their domains.
Besides, $\nabla \psi^*$ and $\nabla \phi^*$ are Lipschitz continuous over $X$ and $Y$.}
\end{assumption}

\begin{remark}
The initial states, including  $x_0 = \nabla \psi^*(u_0)$, $y_0 = \nabla \phi^*(v_0)$, $u_0 = \nabla \psi(x_0)$ and $v_0 = \nabla \phi(y_0)$,  ensure (\ref{alg}) has a unique solution as will be shown in Theorem \ref{the:exist}, and they are satisfied by taking $u_0 \in X$ and $v_0 \in Y$.
{\color{blue} In order to guarantee $(x^\star, y^\star)$ is a solution to (\ref{form}) at the equilibrium point $(u^\star, x^\star, v^\star, y^\star)$ of (\ref{alg}) as will be proved in Lemma \ref{lem:equ}, the extra terms $\nabla \psi(x + \frac t r \dot x)- u$ and $\nabla\phi(y + \frac t r \dot y) - v$ are designed in the evolution of $u$ and $v$.
Additionally, inspired by \cite{antipin1994feedback}, we introduce the derivative feedbacks $\dot x$ and $\dot y$, as damping terms, for the convergence of  (\ref{alg}). }
\end{remark}

\begin{remark}
{\color{blue}
In fact, (\ref{alg}) is a generalization of (\ref{acce:mirror}) in \cite{krichene2015accelerated} under a primal-dual framework.}
Compared with the second-order mirror-descent dynamics in \cite{gao2023second}, (\ref{alg}) is only a first-order algorithm, and hence, $F$ and $G$ are assumed to be strictly convex for its exact convergence.
We focus on a sublinear rate instead of linear convergence, and thus, we do not impose that $F$ and $G$ are strongly convex as \cite{gao2022continuous} did.
While the discounted mirror-descent dynamics only achieved perturbed Nash equilibriums for monotone games \cite{gao2020continuous, gao2022continuous}, (\ref{alg}) reaches the exact saddle point of  (\ref{form}) as will be proved in Theorem \ref{the:con}.
Furthermore, it has a faster rate than the algorithm in \cite{chen2021distributed}.
If the constraints $X$ and $Y$ in (\ref{form}) disappear, (\ref{alg}) degenerates into the primal-dual Nesterov's accelerated algorithms in \cite{zeng2022dynamical, zeng2023distributed} by taking $\nabla \psi^*(u) = u$ and $\nabla \phi^*(v) = v$.
\end{remark}

The introduction of $\psi$ and $\phi$ makes it possible to analyze and design mirror-descent algorithms from the dual space of the Euclidean space.
In practice, they can be employed to deal with some constraints efficiently. We provide the following two examples for illustration, and recommend readers to refer to  \cite{gao2022continuous, gao2023second} for more details.
\begin{enumerate}[(a)]
\item Let $X$ be a nonempty, closed and convex set. Consider
$\psi(x) = \frac 12\Vert x \Vert^2, x \in X$. Then $\nabla \psi^*(u) = P_X (u)$.

\item Let $X$ be a simplex, i.e., $X = \{x = [x_i] \in \mathbb R^p ~|~ x_i \ge 0, \sum_i x_i = 1 \}$. Take $\phi(x) = \sum_i x_i \log(x_i), x \in X$, where $0{\rm log}(0) = 0$. Then
$\nabla \psi^*(u) = \big[\frac {\exp(u_i)}{\sum_i \exp(u_i)}\big]$.
\end{enumerate}

\section{MAIN RESULTS}

This section analyzes properties of (\ref{alg}), including the existence of a solution, the optimality and the convergence.


For convenience, we define $z=[x^{\rm T}, y^{\rm T}]^{\rm T}, w =[u^{\rm T}, 
v^{\rm T}]^{\rm T}$, 
$\nabla R(z)=[\nabla F^{\rm T}(x), \nabla G^{\rm T}(y)]^{\rm T}$,
$\tilde H = [0_{p \times p}, H; -H^{\rm T}, 0_{q \times q}]$,
$\nabla \pi^*(w) = [\nabla \psi^{*{\rm T}}(u), \nabla\phi^{*{\rm T}}(v)]^{\rm T}$, and $\nabla \pi(z) = [\nabla \psi^{\rm T}(x), \nabla\phi^{\rm T}(y)]^{\rm T}$. 
Then dynamics (\ref{alg}) reads as
\begin{equation}
\begin{aligned}
\label{alg:asim}
\dot z =& \frac {r}{t} \big[\nabla \pi^*(w) - z\big] \\
\dot w =& \frac tr \big[\!-\nabla R(z) - \tilde H(z+ \frac tr \dot z) 
+ \nabla \pi(z + \frac t r \dot z) - w\big]
\end{aligned}
\end{equation}
where $z(0) = z_0 = \nabla \pi^*(w_0)$, and $w(0) = w_0 = \nabla \pi(z_0)$.

{\color{blue}Classical results in ordinary differential equations  (ODEs) do not directly imply the existence and uniqueness of a solution to (\ref{alg:asim}) since $r/ t$ is singular at $t = 0$.}
For the analysis, we define an approximate smooth dynamics as
\begin{equation}
\begin{aligned}
\label{alg:apsim}
\dot z =& \frac {r}{\max(\delta,t)} \big[\nabla \pi^*(w) - z\big] \\
\dot w =& \frac tr \big[\!-\nabla R(z) - \tilde H(z+ \frac tr \dot z) 
+ \nabla \pi(z + \frac t r \dot z) - w\big]
\end{aligned}
\end{equation}
where $\delta > 0$, $z(0) = z_0 = \nabla \pi^*(w_0)$, and $w(0) = w_0 = \nabla \pi(z_0)$. Under Assumptions \ref{ass:convex} and \ref{ass:genlip},  (\ref{alg:apsim}) has a unique solution $\big(z_{\delta}(t), w_{\delta}(t)\big)$ by \cite[Th. 3.2]{khalil2002nonlinear}.
In addition, we have the following result, whose proof is given in the Appendix.
\begin{lemma}
\label{lem:equi}
Consider dynamics (\ref{alg:apsim}). 
If Assumptions \ref{ass:convex} and \ref{ass:genlip} hold, then there exists a small enough $t_m$ such that the family of solutions 
$\big((z_{\delta}(t), w_{\delta}(t))|_{[0, t_m]}\big)_{\delta \le t_m}$ 
is equi-continuous and uniformly bounded.
\end{lemma}

First, we establish that there is a unique solution to (\ref{alg:asim}).
\begin{theorem}
\label{the:exist}
If Assumptions \ref{ass:convex} and \ref{ass:genlip} hold,  then dynamics (\ref{alg:asim}) has a unique solution $(z(t), w(t))$ for all $t \ge 0$.
\end{theorem}

\emph{Proof:}
Here we show that there exists a solution to (\ref{alg:asim}). Uniqueness of the solution is proved in the Appendix.

Let $\delta_i = 2^{-i}t_m$ for $i \in \mathbb N$, and $(z_{\delta_i}, w_{\delta_i})$ be a family of solutions to (\ref{alg:apsim}) defined
on $[0, t_m)$.
By Lemma \ref{lem:equi}, the family is equi-continuous and uniformly bounded.
Resorting to the Arzela-Ascoli theorem \cite[Th. 10.1]{royden2010real}, there must be a subsequence
$(z_{\delta_i}, w_{\delta_i})_{i \in \mathcal I}$
that converges uniformly to a limit $(\bar z, \bar w)$, where $\mathcal I \subset \mathbb N$ is an infinite set of indices.

We now prove that $(\bar z, \bar w)$ is a solution to (\ref{alg:asim}).
Clearly, the initial conditions hold.
Take $t_3 \in [0, t_m)$, and let $(\tilde z, \tilde w)$ be a solution to (\ref{alg:asim}) on $t \ge t_3$ with initial condition 
$\big(\bar z(t_3), \bar w(t_3)\big)$.
Since $\lim_{i \to \infty}\big(z_{\delta_i}(t_3), w_{\delta_i}(t_3) \big)_{i\in \mathcal I} = \big(\bar z(t_3), \bar w(t_3)\big)$,
it follows from \cite[Th. 2.8]{teschl2012ordinary} that $(z_{\delta_i}, w_{\delta_i})$ converges to $(\tilde z, \tilde w)$ uniformly on $[t_3, t_3 + \epsilon)$ for some $\epsilon > 0$.
Note that $(z_{\delta_i}, w_{\delta_i})_{i\in \mathcal I}$ also converges to $(\bar z, \bar w)$. Thus,
$(\bar z, \bar w)$ coincides with $(\tilde z, \tilde w)$ on $[t_3, t_3 + \epsilon)$,
and $(\bar z, \bar w)$  is a solution to (\ref{alg:asim}) on $[t_3, t_3 + \epsilon]$.
Due to the arbitrariness of $t_3$, (\ref{alg:asim}) has a solution.
$\hfill\square$


The next lemma addresses the relationship between the saddle point of (\ref{form}) and an equilibrium point of (\ref{alg}).
\begin{lemma}
\label{lem:equ}
Under Assumptions \ref{ass:convex} and \ref{ass:genlip}, $(x^\star, y^\star)$
is the solution to (\ref{form}) if and only if there exists $(u^\star, v^\star)$ such that
$(x^\star, u^\star, y^\star, v^\star)$ is an equilibrium point of dynamics (\ref{alg}).
\end{lemma}

\emph{Proof:}
If $(x^\star, u^\star, y^\star, v^\star)$ {\color{blue}is an equilibrium point} of (\ref{alg}), then
\begin{equation*}
\begin{aligned}
0_p =&~ \nabla \psi^*(u^\star) - x^\star, ~
0_q = \nabla \phi^*(v^\star) - y^\star \\
0_p =& - \nabla F(x^\star) - Hy^\star + 
\nabla \psi(x^\star) - u^\star \\
0_q =& - \nabla G(y^\star) + H^{\rm T} x^\star + 
\nabla \psi(y^\star) - v^\star. 
\end{aligned}	
\end{equation*}

Recalling (\ref{conjugate}) gives $u^\star - \nabla \psi(x^\star) \in \mathcal N_X(x^\star)$ and
$v^\star - \nabla \psi(y^\star) \in \mathcal N_Y(y^\star)$.
Therefore, $- \nabla F(x^\star) - Hy^\star \in \mathcal N_X(x^\star)$, and 
$ - \nabla G(y^\star) + H^{\rm T} x^\star \in \mathcal N_Y(y^\star)$.
It follows from the Karush–Kuhn–Tucker (KKT) optimality conditions \cite[Th. 3.34]{ruszczynski2011nonlinear} that $(x^\star, y^\star)$ is the saddle point of $L$.

Conversely, suppose that $(\tilde x, \tilde y)$ is the saddle point of $L$. Take
$\tilde u - \nabla \psi(\tilde x) = - \nabla F(\tilde x) 
- H \tilde y \in \mathcal N_X(\tilde x)$, and
$\tilde v - \nabla \phi(\tilde y) = -\nabla G(\tilde y) + H^{\rm T} \tilde x
\in \mathcal N_Y(\tilde y)$.
Then $\tilde x = \nabla \psi^*(\tilde u)$, and $\tilde y = \nabla \phi^*(\tilde v)$.
Thus, $(\tilde x, \tilde u, \tilde y, \tilde v)$ is the equilibrium point of (\ref{alg}), and the proof is completed.
$\hfill\square$

Before proving the convergence of (\ref{alg}), we provide a supporting lemma as follows.
\begin{lemma}
\label{lem:set}
Consider dynamics (\ref{alg}). Under Assumptions \ref{ass:convex} and {\ref{ass:genlip}}, $x(t) \in X$, and $y(t) \in Y$ for all $t \ge 0$.
\end{lemma}
\emph{Proof:}
Define $h(t) = \frac 12\Vert x - P_{X}(x)\Vert^2$.
Clearly,
\begin{equation*}
\begin{aligned}
&\dot h(t) =  
\langle x - P_{X}(x), \dot x\rangle
= \frac rt\langle x - P_{X}(x), \nabla \psi^*(u)
- x\rangle \\
&= - \frac rt \Vert x - P_{X}(u)\Vert^2 
+  \frac r t\langle x - P_{X}(x), \nabla \psi^*(u) - P_{X}(x)\rangle.
\end{aligned}	
\end{equation*}

By  (\ref{conjugate}), $\nabla \psi^*(u) \in X$, and then,
$\langle x - P_{X}(x), \nabla \psi^*(u) - P_{X}(x)\rangle \le 0$.
It follows that $\dot h(t) \le 0$. 
Since $x_0 \in X$, $h(t) = 0$, i.e., $x(t) \in X, ~\forall t \ge 0$.
Similarly, $y(t) \in Y, ~\forall t \ge 0$.
This completes the proof.
$\hfill\square$

Finally, we analyze the convergence of (\ref{alg}).
\begin{theorem}
\label{the:con}
Suppose Assumptions \ref{ass:convex} and \ref{ass:genlip} hold.
Let $\big(x(t), \\
u(t), y(t), v(t)\big)$ be a trajectory of (\ref{alg}), and 
$(x^\star, u^\star, y^\star, v^\star)$ be an equilibrium point of (\ref{alg}),
where $r \ge 2$, and $(x^\star, y^\star)$ is the saddle point of  (\ref{form}).
\begin{enumerate}[(i)]
\item The duality gap, defined by $L(x(t), y^\star) - L(x^\star, y(t))$, converges to $0$ with a rate of $\mathcal O(1/t^2)$.

\item The trajectory of $\big(x(t), y(t)\big)$ approaches $(x^\star, y^\star)$.
\end{enumerate}
\end{theorem}
\emph{Proof:}
By Lemma \ref{lem:equ}, $(x^\star, y^\star)$ is the saddle point of  (\ref{form}) for any equilibrium point $(x^\star, u^\star, y^\star, v^\star)$ of (\ref{alg}).
Define 
\begin{equation}
\begin{aligned}
\label{lya}
V =  \frac {t^2}{r}\big(L(x, y^\star) - L(x^\star, y)\big) 
+ rD_{\psi^*}(u, u^\star) + rD_{\phi^\star}(v, v^\star)
\end{aligned}
\end{equation}
where $D_{\psi^*}$ and $D_{\phi^*}$ are the Bregman divergence of $\psi^*$ and $\phi^*$. 
Recalling Lemma \ref{lem:set} gives $L(x, y^\star) - L(x^\star, y) \ge 0$, where the equality holds only if $x=x^\star$ and $y = y^\star$.
Due to the convexity of $\psi^*$ and $\phi^*$,  
$D_{\psi^*}(u, u^\star) + D_{\phi^\star}(v, v^\star) \ge 0$.
Hence, $V \ge 0$, and moreover, $V = 0$ only if $(x, y) = (x^\star, y^\star)$.

Combining (\ref{alg}) with (\ref{lya}), we obtain
\begin{equation*}
\begin{aligned}
\dot V &= \frac {2t}{r}\big(L(x, y^\star) - L(x^\star, y)\big) 
+ t \langle \nabla F(x) + H y^\star, \\
&~~~ \nabla \psi^*(u) - x \rangle 
+ t \langle \nabla G(y) - H^{\rm T} x^\star, \nabla \phi^*(v) - y \rangle \\
& + r \langle \nabla\psi^*(u) - \nabla \psi^*(u^\star), \dot u \rangle 
+ r \langle \nabla\phi^*(v) - \nabla \phi^*(v^\star), \dot v\rangle.
\end{aligned}
\end{equation*}
Note that
\begin{equation*}
\begin{aligned}
t \langle \nabla F(x) + H y^\star, \nabla \psi^*(u) - x \rangle 
+ r \langle \nabla\psi^*(u) -\nabla\psi^*(u^\star), \dot u \rangle & \\
=   t \langle x + \frac tr \dot x - x^\star, - H (y + \frac tr \dot y - y^\star) + \nabla \psi(x + \frac t r \dot x) - u \rangle & \\
- t\langle  x - x^\star, \nabla F(x) + H y^\star\rangle &
\end{aligned}
\end{equation*}
and moreover,
\begin{equation*}
\begin{aligned}
t \langle \nabla G(y) - H^{\rm T} x^\star, \nabla \phi^*(v) - y \rangle 
+ r \langle \nabla\phi^*(v) \!-\! \nabla\phi^*(v^\star), \dot v \rangle & \\
= t \langle y + \frac tr \dot y - y^\star, 
 H^{\rm T} (x + \frac tr \dot x - x^\star) 
+ \nabla \phi(y + \frac t r \dot y) - v \rangle & \\
- t \langle y - y^\star,  \nabla G(y) - H^{\rm T} x^\star\rangle.&
\end{aligned}
\end{equation*}
Consequently, 
\begin{equation}
\begin{aligned}
\label{pf:eq0}
\dot V & = \frac {2t}{r}\big(L(x, y^\star) - L(x^\star, y)\big) \\
& -t \langle  x - x^\star, \nabla F(x) + H y^\star \rangle  - t \langle y - y^\star, \nabla G(y) - H^{\rm T} x^\star \rangle \\
& + t\langle \nabla\psi^*(u) - x^\star, 
\nabla \psi(x + \frac t r \dot x) - u \rangle \\
& + t \langle \nabla\phi^*(v) - y^\star, 
\nabla \phi(y + \frac t r \dot y ) - v \rangle.
\end{aligned}
\end{equation}

Since $x +  \frac t r \dot x = \nabla \psi^*(u)$,  
$\nabla \psi(x + \frac t r \dot x) - u \in - \mathcal N_X(x + \frac t r \dot x)$ by  Lemma \ref{lem:gen}.
Therefore,
\begin{equation*}
\begin{aligned}
 \langle &\nabla\psi^*(u) - x^\star, 
\nabla \psi(x +  \frac t r\dot x) - u \rangle  \\
&= \langle x +  \frac t r \dot x - x^\star, 
\nabla \psi(x +  \frac t r \dot x) - u \rangle  \le 0.
\end{aligned}
\end{equation*}

Similarly, 
$\langle \nabla\phi^*(v) - y^\star, 
\nabla \phi(y + \frac t r \dot y) - v \rangle \le 0.$
Due to the convexity of $F$ and $G$, we obtain
\begin{equation*}
\begin{aligned}
L(x^\star, y) - L(x, y^\star) \ge  \langle  x^\star - x, \nabla F(x) + H y^\star\rangle &\\ 
+ \langle y^\star - y, \nabla G(y) - H^{\rm T} x^\star\rangle&.
\end{aligned}
\end{equation*}

By (\ref{pf:eq0}), $\dot V \le  (\frac {2t}{r} - t)\big(L(x, y^\star) - L(x^\star, y)\big)$. 
As a result, $\dot V \le 0$ if $r \ge 2$.
Let $m_0 = V(x_0, u_0, y_0, v_0)$.
Then $V \le m_0, \forall t \ge 0$. 
Furthermore,
$L(x, y^\star) - L(x^\star, y) \le \frac {r}{t^2} V  \le \frac {r}{t^2} m_0$,
i.e.,
$L (x, y^\star)  - L (x^\star, y)$ converges to $0$ with a rate of $\mathcal O(1/t^2)$.
Therefore, part (i) is proved.

It is clear that
$L(x, y^\star) - L(x^\star, y^\star) \le L(x, y^\star) - L(x^\star, y) \le \frac {r}{t^2} m_0$.
Since $F$ is strictly convex, $L(x, y^\star)$ has a unique minimizer, and as a result, $x(t)$ converges to $x^\star$.
By a similar procedure, $y(t)$ approaches $y^\star$.
Thus, part (ii) is proved, and the proof is completed.
$\hfill\square$

Theorem \ref{the:con} indicates that (\ref{alg}) reaches the exact saddle point of (\ref{form}) with a rate of $\mathcal O(1/t^2)$.
{\color{blue}
It has a faster performance than the primal-dual methods in \cite{
zeng2018distributed, lei2016primal, gharesifard2013distributed, yang2016multi}.
Besides, as a first-order algorithm, 
(\ref{alg}) has an optimal rate in the worst case \cite[Th.  2.1.7]{nesterov2003introductory}.}

{\color{blue}
\begin{remark}
In fact, the order of the convergence rate of (\ref{alg}) remains unchanged if its right-hand side is multiplied by a positive constant gain. However, time-varying gains should be carefully chosen because they may produce unbounded variable derivatives and make algorithms impractical.
\end{remark}
}

{\color{blue}
\begin{remark}
The continuous-time algorithm (\ref{alg}) may be discretized for its numerical implementation.
Referring to \cite{shi2019acceleration}, there are various discretization schemes for an ODE. By numerical tests, we have observed that 
using a similar discretization scheme as that of \cite{krichene2015accelerated}, the corresponding discrete-time algorithm from (\ref{alg}) still maintains the acceleration characteristics. The convergence, robustness and numerical stability of the scheme will be investigated in the future.
\end{remark}
}

\section{DISTRIBUTED APPLICATIONS}

In this section, we apply the proposed algorithm (\ref{alg}) to constrained distributed optimization and two-network bilinear zero-sum games.

\subsection{Constrained Distributed Optimization}

Consider a network of $n$ agents interacting over an undirected and connected graph $\mathcal G(\mathcal V, \mathcal E)$, where agent $i$ only knows a local function $f_i: \Omega_i \rightarrow \mathbb R$ and a closed set $\Omega_i \subset \mathbb R^p$, and all agents cooperate to
find a consensus solution that minimizes the global function $\sum_{i \in \mathcal V} f_i(x)$. 
The problem is formulated as
\begin{equation}
\begin{aligned}
\label{Ex1}
\min ~ f(\bm x) := \sum\nolimits_{i\in \mathcal V} f_i(x_i) ~~
{\rm s.t.}~\bm L \bm x = 0, ~ \bm x \in \Omega
\end{aligned}
\end{equation}
where $\bm x = [x_1^{\rm T}, \dots, x_n^{\rm T}]^{\rm T} \in \mathbb R^{pn}$,
$\mathcal L \in \mathbb R^{n \times n}$ is the Laplacian matrix of $\mathcal G$,
$\bm L = \mathcal L \otimes I_p \in \mathbb R^{pn \times pn}$ and $\Omega = \Omega_1 \times \dots \times \Omega_n$.

We make the following assumption for (\ref{Ex1}).
\begin{assumption}
\label{ass:ex1}
For $i \in \mathcal V$, $\Omega_i$ is closed and convex,
$f_i$ is continuously differentiable and {\color{blue}strictly convex} on an open set containing $\Omega_i$, and $\nabla f_i$ is Lipschitz continuous over $\Omega_i$.
{\color{blue}The Slater's condition holds}, and there exists a unique solution.
\end{assumption}

{\color{blue}
To solve (\ref{Ex1}), distributed algorithms have been explored in \cite{lei2016primal, nedic2010constrained}} with applications in
motion planning, alignment of multiple vehicles and distributed optimal consensus of multi-agent systems.
Due to the presence of $\Omega$, the accelerated protocol in \cite{zeng2022dynamical} cannot be adopted.

Referring to \cite{lei2016primal}, we first reformulate (\ref{Ex1}) as
\begin{equation}
\begin{aligned}
\label{ex1:Larg}
\min_{\bm x \in \Omega}~\max_{\bm y \in \mathbb R^{pn}}~
S_1 (\bm x, \bm y) := f(\bm x) 
+ \bm y^{\rm T} \bm L \bm x + \frac 12 \bm x^{\rm T} \bm L \bm x
\end{aligned}
\end{equation}
where $\bm y = [y_1^{\rm T}, \dots, y_n^{\rm T}]^{\rm T} \in \mathbb R^{pn}$ is the dual variable. 
Then we introduce variables $\bm u = [u_1^{\rm T}, \dots, u_n^{\rm T}]^{\rm T} \in \mathbb R^{pn}$ and $\bm v = [v_1^{\rm T}, \dots, v_n^{\rm T}]^{\rm T} \in \mathbb R^{pn}$, and generating functions $\psi_i$ with  $\nabla \psi_i^*(u_i) = P_{\Omega_i}(u_i)$. 
Following that, the accelerated distributed algorithm {\color{blue} is given by} 
\begin{equation}
\begin{aligned}
\label{alg:ex1}
\dot x_i =& \frac rt \big[P_{\Omega_i}(u_i) - x_i  \big],~
\dot y_i = \frac rt \big[v_i - y_i \big]  \\
\dot u_i =& \frac tr \big[\!-\nabla_{x_i} S_1(\bm x, \bm y + \frac t r \dot {\bm y})
+ x_i + \frac tr \dot x_i - u_i\big] \\
\dot v_i =& \frac tr \big[\nabla_{y_i} S_1(\bm x 
+ \frac t r \dot {\bm x}, \bm y) + y_i + \frac t r \dot y_i - v_i\big]
\end{aligned}	
\end{equation}
where 
$x_i(0) = x_{i, 0} = u_{i, 0}$,
$u_{i}(0) = u_{i, 0} \in \Omega_i$, 
$y_i(0) = y_{i, 0} = v_{i, 0}$, and
$v_i(0) = v_{i, 0} \in \mathbb R^p$.

Finally, we establish the convergence of (\ref{alg:ex1}) as follows.
\begin{corollary}
\label{Cor:Ex1}
Suppose  Assumption \ref{ass:ex1} holds. Let $\big(\bm x(t),  \bm u(t),\\
\bm y(t), \bm v(t) \big)$ be a trajectory of (\ref{alg:ex1}), and $(\bm x^\star ,  \bm u^\star, \bm y^\star, \bm v^\star)$ be an equilibrium point  of (\ref{alg:ex1}), where $r \ge 2$, $(\bm x^\star, \bm y^\star)$ is a saddle point of $S_1$,
and $\bm x^\star$ is the optimal solution to (\ref{Ex1}).
\begin{enumerate}[(i)]
\item The duality gap, defined by $S_1 \big(\bm x(t), \bm y^\star \big) - S_1 \big(\bm x^\star, \bm y(t) \big)$, converges to $0$ with a rate of $\mathcal O(1/t^2)$.

\item The trajectory of {\color{blue}$\bm x(t)$ approaches $\bm x^\star$.}

\item {\color{blue} Both $\Vert \bm L \bm x \Vert$ and $|f(\bm x) - f(\bm x^\star)|$ converge to $0$ with a rate of $\mathcal O(1/t)$.}
\end{enumerate}
\end{corollary}
\emph{Proof:}
By a similar procedure as the proof of Theorem \ref{the:con},
items (i) and (ii) hold resorting to a function defined as
$$
V_1 = \frac{t^2}{r} \big(S_1(\bm x, \bm y^\star) - S_1(\bm x^\star, \bm y) \big)
+ rD_{\psi^*}(\bm u, \bm u^\star) + r D_{\phi^*}(\bm v, \bm v^\star)
$$
where $D_{\psi^*}(\bm u, \bm u^\star) = \sum_{i \in \mathcal V} D_{\psi_i^*}(u_i, u_i^\star)$, 
$D_{\phi^*}(\bm v, \bm v^\star) = \sum_{i \in \mathcal V} D_{\phi_i^*}(v_i, v_i^\star)$,
$D_{\psi_i^*}$ and $D_{\phi_i^*}$ are the Bregman divergence, $\psi_i(x_i) = \frac 12 \Vert x_i\Vert^2, x_i \in \Omega_i$, and $\phi_i(y_i) = \frac 1 2 \Vert y_i \Vert^2, y_i \in \mathbb R^p$.

{\color{blue}
Since $S_1(\bm x, \bm y^\star ) - S_1(\bm x^\star , \bm y) 
\ge \frac 12 \bm x^T \bm L \bm x$, it follows from (i) that
$\Vert \bm L \bm x \Vert$ approaches $0$ with a rate of $\mathcal O(1/t)$. In addition, $0 \le f(\bm x) + \bm y^{\star T} \bm L \bm x - f(\bm x^*) \le S_1(\bm x, \bm y^*) - S_1(\bm x^*, \bm y)$, and as a result,
$$
|f(\bm x) - f(\bm x^\star )| \le \Vert \bm y^\star \Vert \cdot \Vert \bm L \bm x\Vert + S_1(\bm x, \bm y^\star ) - S_1(\bm x^\star , \bm y).
$$
Therefore, $|f(\bm x) - f(\bm x^\star )|$ converges to $0$ with a rate of $\mathcal O(1/t)$, and item (iii) holds.
The proof is completed}.
$\hfill\square$

{\color{blue}
\begin{remark}
We should mention that (\ref{Ex1}) has a unique solution $\bm x^\star$, but admits multiple dual variables $\bm y^*$ such that $(\bm x^\star, \bm y^\star)$ is a saddle point of $S_1$ \cite{hamedani2021primal}.
Thus, we only analyze the converge of $\bm x(t)$.
Referring to \cite{xu2017accelerated, luo2024universal}, 
we convert the convergence result on the duality gap into the constraint violation and the primal suboptimality in Corollary \ref{Cor:Ex1} (iii).
\end{remark}
}

\subsection{Two-Network Bilinear Zero-Sum Games}

Consider two undirected and connected networks {\color{blue}$\mathcal G_1(\mathcal V_1, \mathcal E_1)$ and $\mathcal G_2(\mathcal V_2, \mathcal E_2)$ with $n_1$ and $n_2$ agents}, where agent $i$ in $\mathcal G_1$ knows a function $f_i: X_i \rightarrow \mathbb R$ and a constraint set $X_i \subset \mathbb R^p$, and moreover, agent $j$ in $\mathcal G_2$ knows a function $g_j: Y_j \rightarrow \mathbb R$ and a constraint set $Y_j \subset \mathbb R^q$.
Let $X = X_1 \times \cdots \times X_{n_1}$ and
$Y = Y_1 \times \cdots \times Y_{n_2}$.
The payoff function $U(\bm x, \bm y): X
\times Y \rightarrow \mathbb R$ of the zero-sum game is
$$
U(\bm x, \bm y) := 
\tilde f(\bm x) + \bm x^{\rm T} B \bm y - \tilde g(\bm y)
$$
where $\bm x = [x_1^{\rm T}, \dots, x^{\rm T}_{n_1}]^{\rm T}$, 
$\bm y = [y^{\rm T}_1, \dots, y^{\rm T}_{n_2}]^{\rm T}$, 
$\tilde f(\bm x) = \sum_{i \in \mathcal V_1} f_i(x_i)$,
$\tilde g(\bm y) = \sum_{j \in \mathcal V_2} g_j(y_j)$,  
$\bm x^{\rm T} B \bm y = \sum_{i \in \mathcal V_1}\sum_{j \in \mathcal V_2} x_i^{\rm T} H_{ij} y_j$, and $H_{ij} \in \mathbb R^{p \times q}$ is a non-zero matrix only if agent $i$ of $\mathcal G_1$ and  agent $j$ of $\mathcal G_2$ can observe each other's decision variable. 
Feasibility sets of $\mathcal G_1$ and $\mathcal G_2$ are
$\Omega_1 = \{\bm x \in X: x_1 = \dots = x_{n_1}\}$ and
$\Omega_2 = \{\bm y \in Y: y_1 = \dots = y_{n_2}\}$.
Seeking {\color{blue}an} Nash equilibrium of the game can be cast into
\begin{equation}
\begin{aligned}
\label{Ex2}
\min_{\bm x \in X} \max_{\bm y \in Y}~ U(\bm x, \bm y)
~~{\rm s.t.}~ \bm L_1 \bm x = 0_{pn_1},~ 
\bm L_2 \bm y = 0_{qn_2}
\end{aligned}
\end{equation}
where $\mathcal L_1$ and $\mathcal L_2$ are the Laplacian matrices of $\mathcal G_1$ and $\mathcal G_2$, $\bm L_1 = \mathcal L_1 \otimes I_p \in \mathbb R^{pn_1 \times pn_1}$, and $\bm L_2 = \mathcal L_2 \otimes I_q \in \mathbb R^{qn_2 \times qn_2}$.

We make the following assumption for (\ref{Ex2}).
\begin{assumption}
\label{ass:ex2}
For each $i \in \mathcal V_1$ and $j \in \mathcal V_2$, $f_i$ and $g_j$ are continuously differentiable and strictly convex on two open sets containing $X_i$ and $Y_j$, respectively.
Their gradients are Lipschitz continuous over $X_i$ and $Y_j$. 
Moreover, {\color{blue}the Slater's condition holds}, and there exists a unique Nash equilibrium.
\end{assumption}

Problem (\ref{Ex2}) models distributed optimization with parameter uncertainties, and  arises in distributed adversarial resource allocation of multiple communication channels \cite{gharesifard2013distributed}.

It follows from \cite{zeng2023distributed} that (\ref{Ex2}) can be reformulated as
\begin{equation}
\begin{aligned}
\min_{\bm x, \bm \mu} ~
\max_{\bm y, \bm \lambda}~
S_2(\bm x, \bm \lambda, \bm y, \bm \mu)
: = U(\bm x, \bm y) + \bm \lambda^{\rm T} \bm L_1 \bm x& \\
- \bm \mu^{\rm T} \bm L_2 \bm y + \frac 12 \bm x^{\rm T} \bm L_1 \bm x
- \frac 12 \bm y^{\rm T} \bm L_2 \bm y&
\end{aligned}
\end{equation}
where $\bm \lambda = [\lambda_1^{\rm T}, \dots, \lambda_{n_1}^{\rm T}]^{\rm T} \in \mathbb R^{p n_1}$ and $\bm \mu = [\mu_1^{\rm T}, \dots, \mu_{n_2}^{\rm T}]^{\rm T} \in \mathbb R^{q n_2}$ are dual variables for the consensus constraints.

We introduce variables 
$\bm u = [u_1^{\rm T}, \dots, u_{n_1}^{\rm T}]^{\rm T} \in \mathbb R^{p n_1}$,
$\bm \tau = [\tau_1^{\rm T}, \dots, \tau_{n_1}^{\rm T}]^{\rm T} \in \mathbb R^{p n_1}$,
$\bm v = [v_1^{\rm T}, \dots, v_{n_2}^{\rm T}]^{\rm T} \in \mathbb R^{q n_2}$,
and $\bm \nu = [\nu_1^{\rm T}, \dots, \nu_{n_2}^{\rm T}]^{\rm T} \in \mathbb R^{q n_2}$.
{\color{blue}Consider the corresponding distributed accelerated algorithm from (\ref{alg}) as}
\begin{equation}
\begin{aligned}
\label{alg:ex2}
\dot x_i =& \frac rt \big[\nabla \psi_i^*(u_i) - x_i \big], ~
\dot \lambda_i = \frac rt \big[\tau_i - \lambda_i \big]  \\
\dot u_i =& \frac tr \big[\! -\! \nabla_{x_i}\! S_2\big(\bm x, \bm \lambda \! + \! \frac t r \dot{\bm \lambda}, y \! + \! \frac t r \dot{\bm y}, \bm \mu \big) \! + \! \nabla \psi_i \big(x_i \! + \! \frac t r \dot x_i\big) \! -\! u_i\big] \\
\dot \tau_i =& \frac tr \big[\nabla_{\lambda_i} S_2\big(\bm x + \frac tr \dot{\bm x}, \bm \lambda, \bm y, \bm \mu \big) + \lambda_i + \frac t r\dot{\lambda_i} - \tau_i \big] \\
\dot y_j =& \frac rt \big[\nabla \phi_j^*(v_j) - y_j \big], ~
\dot \mu_j = \frac rt \big[\nu_j - \mu_j \big] \\
\dot v_j =& \frac tr \big[\nabla_{y_j} \! S_2 \big(\bm x \! + \! \frac t r \dot{\bm x},  \bm \lambda,  \bm y, \bm \mu \! + \! \frac t r \dot{\bm \mu} \big) \! + \! \nabla \phi_j(y_j \! + \! \frac t r \dot y_j) \! - \! v_j\big] \\
\dot \nu_j =& \frac tr \big[\! -\nabla_{\mu_i} S_2 \big(\bm x, \bm \lambda, \bm y  + \frac tr \dot{\bm y}, \bm \mu \big) + \mu_j +\frac tr \dot{\mu_j} - \nu_j\big]
\end{aligned}	
\end{equation}
where $x_i(0) \!=\! x_{i, 0} \!=\! \nabla \psi_i^*(u_{i, 0})$, $u_i(0) \!=\! u_{i, 0} \in X_i$,
$\lambda_i(0) \!=\! \lambda_{i, 0} \!=\! \tau_{i, 0}$, $\tau_i(0) \!=\! \tau_{i, 0} \in \mathbb R^p$,
$y_j(0) \!=\! y_{j, 0} \!=\! \nabla \phi_i^*(v_{j, 0})$, 
$v_j(0) \!=\! v_{j, 0} \in Y_j$,
$\mu_j(0) \!=\! \mu_{j, 0} \!=\! \nu_{i, 0}$ and $\nu_j(0) \!=\! \nu_{j, 0} \in \mathbb R^q$. The generating functions $\psi_i$ and $\phi_j$, depending on $X_i$ and $Y_j$,  are manually selected with similar requirements as that of $\psi$ and $\phi$ in Assumption \ref{ass:genlip}.
Clearly, (\ref{alg:ex2}) covers the algorithm in \cite{zeng2023distributed}, but $\psi_i$ and $\phi_i$ are introduced to handle constraints.
Its convergence is addressed as follows.
\begin{corollary}
Suppose Assumption \ref{ass:ex2} holds.
Let $\big(\bm x(t), \bm u(t), \\
\bm \lambda(t), \bm \tau(t), \bm y(t), \bm v(t), \bm \mu(t), \bm \nu(t) \big)$ be a trajectory of (\ref{alg:ex2}), and $(\bm x^\star, \bm u^\star, \bm \lambda^\star, \bm \tau^\star, \bm y^\star, \bm v^\star, \bm \mu^\star, \bm \nu^\star)$ be an equilibrium point of (\ref{alg:ex2}).
where $r \ge 2$, $(\bm x^\star, \bm \lambda^\star, \bm y^\star, \bm \mu^\star)$ is a saddle point of $S_2$, and 
$(\bm x^\star, \bm y^\star)$ is the NE of (\ref{Ex2}).
\begin{enumerate}[(i)]
\item The duality gap, defined by $S_2 \big(\bm x(t), \bm \lambda^\star, \bm y^\star, \bm \mu (t)\big) - S_2(\bm x^\star, \bm \lambda(t), \bm y(t), \bm \mu^\star)$, converges to $0$ with a rate of $\mathcal O(1/t^2)$.

\item The trajectory of $\big(\bm x(t), \bm y(t)\big)$ approaches $(\bm x^\star, \bm y^\star)$.

\item {\color{blue}All $\Vert \bm L_1 \bm x \Vert$, $\Vert \bm L_2 \bm y \Vert$
and $|U(\bm x, \bm y^\star) - U(\bm x^\star, \bm y)|$ converge to $0$ with a rate of $\mathcal O(1/t)$.}	
\end{enumerate}
\end{corollary}

\emph{Proof:}
The proof of items (i) and (ii) is similar to that of Theorem \ref{the:con} by defining a function as
\begin{equation*}
\begin{aligned}
V_2 =& \frac {t^2}{r}\big(S_2(\bm x, \bm \lambda^\star, \bm y^\star, \bm \mu) - S_2(\bm x^\star, \bm \lambda, \bm y, \bm \mu^\star) \big) + r D_{\psi_{u}^*}(\bm u, \bm u^\star) \\
& + r D_{\psi_{\tau}^*}(\bm \tau, \bm \tau^\star) 
+ r D_{\phi_{v}^*}(\bm v, \bm v^\star) 
+ r D_{\phi_{\nu}^*}(\bm \nu, \bm \nu^\star)
\end{aligned}	
\end{equation*}
where 
$D_{\psi_{u}^*}(\bm u, \bm u^\star) = \sum_{i \in \mathcal V_1} D_{\psi_{i}^*}(u_i, u_i^\star)$,
$D_{\psi_{\tau}^*}(\bm \tau, \bm \tau^\star) = \sum_{i \in \mathcal V_1} D_{\widetilde \psi_i^*}(\tau_i, \tau_i^\star)$
$D_{\phi_{v}^*}(\bm v, \bm v^\star)  = \sum_{j \in \mathcal V_2} 
D_{ \phi_j^*}(v_j, v_j^\star)$,
$D_{\phi_{\nu}^*}(\bm \nu, \bm \nu^\star)  = \sum_{j \in \mathcal V_2} D_{\widetilde \phi_j^*}(\nu_j, \nu_j^\star)$,
$\widetilde \psi_i(\lambda_i) = \frac 12 \Vert\lambda_i\Vert^2, \lambda_i \in \mathbb R^p$, and
$\widetilde \phi_j(\mu_j) = \frac 12 \Vert \mu_j \Vert^2, \mu_j \in \mathbb R^q$.
 The proof of item (iii) can refer to that of Corollary \ref{Cor:Ex1} (iii).  
 $\hfill\square$

\subsection{Numerical Simulations}

Here we provide two examples for illustration.

\begin{example}
{\color{blue}
Consider problem (\ref{Ex1}), where $n = 20$, each node of $\mathcal G$ connects with the others with a probability of $0.7$. For $i \in \mathcal V$, 
$f_i$ is a logistic regression  function given by $f_i(x_i) = \frac {1}{m_i} \sum_{j = 1}^{m_i} \log\big[1 + \exp\big(\! -(a_{ij}^T \hat x_i + x_{i, 0}) l_{ij} \big)\big]$, and $\Omega _i= \{x_i \in \mathbb R^p: \Vert x_i \Vert \le c\}$, where $x_i = [x_{i, 0}, \hat x_i^T]^T \in \mathbb R^p$,
$p = 10$, $m_i = 5$, $c = 100$,  all entries of 
$a_{ij} \in \mathbb R^{p-1}$ are randomly drawn from $[0, 1]$, 
and $l_{ij}$ is a random number in $\{-1, +1\}$.
We solve the problem by the accelerated algorithm (\ref{alg:ex1}) and the primal-dual method in \cite{yang2016multi}, where we approximate  trajectories of the continuous-time dynamics using function ode45 of Matlab.
Fig. 1(a) shows the trajectories of $\log[S_1(\bm x, \bm y^\star) - S_1(\bm x^\star, \bm y)]$, and indicates that (\ref{alg:ex1}) has a faster rate than the standard primal-dual algorithm.
It is also clear that for (\ref{alg:ex1}), the trajectory of the duality gap is not guaranteed to be monotone with respect to time, which is similar to that of Nesterov's acceleration for unconstrained optimization. }
Fig. 1(b) implies that all agents reach a consensus solution due to $\lim_{t \to \infty} \Vert \bm L \bm x \Vert = 0$.
\begin{figure}[htp]
\centering
\includegraphics[scale=0.26]{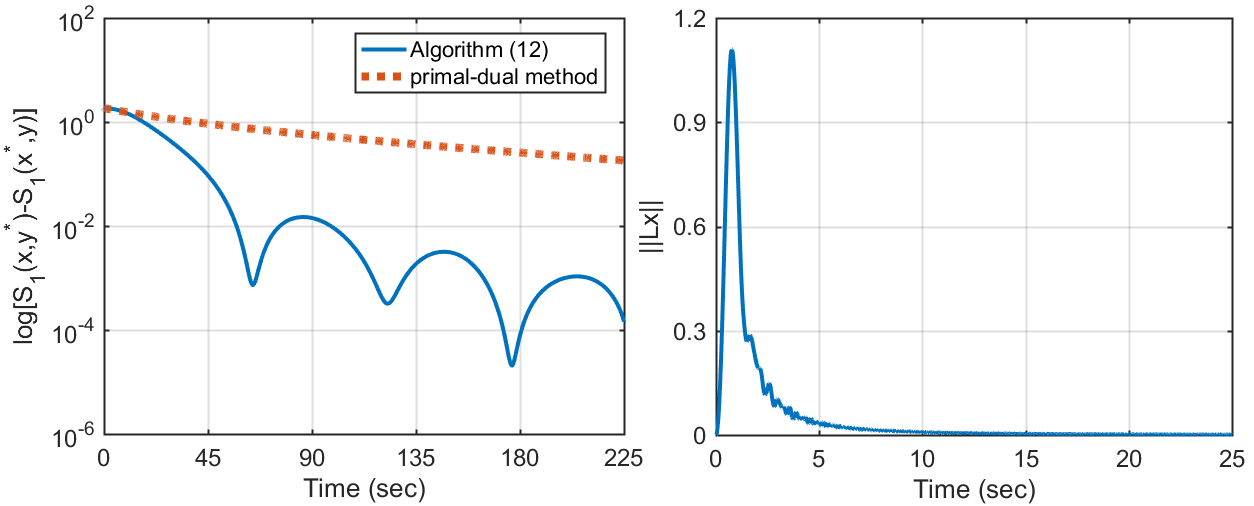}
\caption{(a) Comparative results for  (\ref{alg:ex1}) and the primal-dual method in \cite{yang2016multi}. (b) The trajectory of $\Vert \bm L \bm x \Vert $.}
\label{Fig1}
\end{figure}
\end{example}

\begin{example}
{\color{blue}
Consider problem (\ref{Ex2}), where $n_1 = n_2 = 10$, $p = q = 20$, $B$ is an identity matrix, each node of $\mathcal G_1$ connects with the others with a probability of $0.6$, and $\mathcal G_2$ is generated by a similar procedure as that of $\mathcal G_1$.
Each $f_i$ is a linear function as $f_i(x_i) = a_i^{\rm T} x_i + b_i$,
$g_j$ is a log-sum-exp function as 
$g_j(y_j) =  \rho \log \big[\sum_{k = 1}^m \exp (c_{jk}^T y_j - d_{jk}) / \rho \big]$, 
$X_i = \{x_i  \in \mathbb R^p~|~l_i \le x_i \le u_i\}$, 
and $Y_j = \mathbb R^q$,
where entries of $a_i \in \mathbb R^p$ and $b_i$ are randomly drawn from $[0, 1]$, 
entries of $l_i$ and $u_i$ are random vectors drawn from $[-2, 0]$ and $[1, 5]$, $\rho = 20$,
$m = 8$, $c_{jk} \in \mathbb R^q$ is a random vector from i.i.d standard Guassian distribution,
and $d_{jk} \in \mathbb R$ is a random number from Guassian distribution with mean $0$ and variance $2$.}
Fig. 2(a) presents the trajectories of $U(\bm x, \bm y)$ under dynamics (\ref{alg:ex2}) and the primal-dual method in \cite{gharesifard2013distributed}, and indicates that (\ref{alg:ex2}) converges with a faster rate.
Fig. 2(b) implies that the constraints of (\ref{Ex2}) are satisfied as $t$ tends to infinity.
\begin{figure}[htp]
\centering
\includegraphics[scale=0.26]{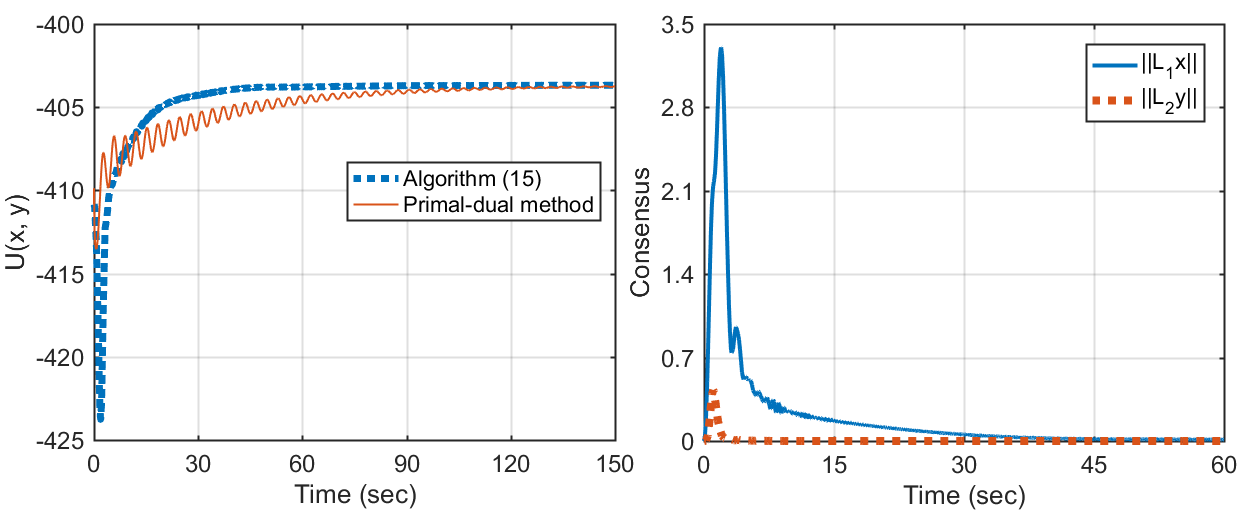}
\caption{(a)  Comparative results for  (\ref{alg:ex2}) and the primal-dual method in \cite{gharesifard2013distributed}. (b) Trajectories of $\Vert \bm L_1 \bm x \Vert$ and $\Vert \bm L_2 \bm y \Vert$.}
\end{figure}
\end{example}

\section{CONCLUSION}

This paper considered solving a class of constrained saddle point problems.
Inspired by a primal-dual framework, we proposed a novel continuous-time accelerated mirror-descent algorithm.
We showed that the algorithm converged with a rate of $\mathcal O(1/t^2)$, which was faster than the standard primal-dual method.
{\color{blue}Following that}, two distributed accelerated algorithms were obtained for constrained distributed optimization and two-network bilinear zero-sum games.
Numerical simulations were also provided for verification.

\appendix

\emph{Proof of Lemma \ref{lem:equi}:}
Under Assumptions \ref{ass:convex} and \ref{ass:genlip}, (\ref{alg:apsim}) has a unique solution on $[0, t_0]$ for some $t_0 > 0$.
Define
$A_\delta(t)= \sup_{0\le s\le t} \Vert \dot w_{\delta}(s) \Vert / s$,
$B_\delta(t)=  \sup_{0\le s\le t} \Vert z_{\delta}(s) - z_0 \Vert / s$,
and $C_\delta(t)= \sup_{0\le s\le t} \Vert \dot z_{\delta}(s) \Vert$. It is clear that 
\begin{equation*}
\begin{aligned}
\label{ine:wd11}
\Vert\dot w_{\delta}(t) \Vert & \le 
\frac tr \Vert \nabla R(z_{\delta}(t))\Vert
+ \frac tr \lambda_H \Vert z_{\delta}(t) \Vert +
\frac {t^2}{r^2} \lambda_H \Vert\dot z_{\delta}(t) \Vert \\
& + \frac t r \Vert \nabla \pi(z_\delta(t) 
+ \frac t r \dot z_\delta(t)) - w_\delta(t) \Vert
\end{aligned}
\end{equation*}
where $\lambda_H$ is the maximal singular value of $\tilde H$. 

Let $L_{\pi}$, $L_{\pi^*}$ and $L_R$ be the Lipschitz constants of $\nabla \pi$, $\nabla \pi^*$ and $\nabla R$ over $X \times Y$.
{\color{blue}
We first show that $z_\delta(t) + \frac{t}{r}\dot z_{\delta}(t) \in X \times Y$ for all $t \ge 0$.
If $t \ge \delta$, $z_\delta(t) + \frac t r \dot z_\delta(t) = \nabla \pi^*(w_\delta(t)) \in X \times Y$ by Lemma 1. 
If $t < \delta$,  $z_\delta(t) + \frac \delta r \dot z_\delta(t) = \nabla \pi^*(w_\delta(t)) \in X \times Y$.
Besides, it holds that $z_\delta(t) \in X \times Y$ for all $t \ge 0$ by a similar procedure as the proof of Lemma \ref{lem:set}.
Due to the convexity of $X \times Y$,
$z_\delta(t) + \frac t r  \dot z_\delta(t) =
\big(1 - \frac t \delta \big) z_{\delta}(t) + \frac t \delta \big(z_{\delta}(t) + \frac \delta r \dot z_\delta(t)\big) \in X \times Y$.
As a result, $z_\delta(t) + \frac t r  \dot z_\delta(t) \in X \times Y,
\forall t \ge 0$.}
Recalling $w_0 = \nabla \pi(z_0)$ yields
\begin{equation*}
\begin{aligned}
\Vert &\nabla \pi(z_\delta(t)  + \frac t r \dot z_\delta(t)) - w_\delta(t)\Vert\\ 
& \le L_{\pi}\Vert z_\delta(t) - z_0 \Vert 
+ \frac t r L_{\pi} \Vert \dot z_\delta(t) \Vert 
+ \Vert w_\delta(t) - w_0 \Vert.
\end{aligned}
\end{equation*}
Then
\begin{equation*}
\begin{aligned}
\label{ine:wd1}
& \Vert  \dot w_{\delta}(t) \Vert  \le 
\frac tr L_R \Vert z_{\delta}(t) - z_0\Vert 
+ \frac tr\Vert R(z_0)\Vert \\
&  + \frac tr \lambda_H \Vert z_{\delta}(t) 
- z_0 \Vert 
+ \frac tr \lambda_H \Vert z_0 \Vert 
+ \frac {t^2}{r^2} \lambda_H \Vert\dot z_{\delta}(t) \Vert\\
& + \frac {t^2}{r^2} L_{\pi} \Vert \dot z_\delta(t) \Vert 
+ \frac t r L_{\pi}\Vert z_\delta(t) - z_0 \Vert 
+ \frac t r \Vert w_\delta(t) - w_0 \Vert.
\end{aligned}
\end{equation*}

Since $\Vert w_\delta(t) - w_0 \Vert \le t^2 A_{\delta} (t) /2$,
\begin{equation}
\begin{aligned}
\label{A_ineq}
r A_{\delta}(t) \le \beta_0 + \beta_1 t B_\delta(t) + \beta_2 t C_{\delta}(t)
+ \beta_3 t^2 A_{\delta}(t)
\end{aligned}
\end{equation}
where $\beta_0= \Vert R(z_0) \Vert + \lambda_H \Vert z_0 \Vert$,  $\beta_1 = L_R + \lambda_H + L_{\pi}$, 
$\beta_2 = (L_{\pi} +\lambda_H)/r$
and $\beta_3 = 1/2$.

It follows from (\ref{alg:apsim}) that
if $t \le \delta$,
\begin{equation*}
\begin{aligned}
e^{\frac {rt}{\delta}} \big( 
\dot z_{\delta} +\frac r\delta(z_{\delta} - z_0)\big)=
\frac r \delta e^{\frac {rt}{\delta}} \big(\nabla\pi^*(w_\delta)
- \nabla \pi^*(w_0)\big)
\end{aligned}
\end{equation*}
i.e.,
$\frac {\rm d}{\rm dt} e^{\frac {rt}{\delta}}\big( z_{\delta} - z_0\big) =
\frac r \delta e^{\frac {rt}{\delta}} \big(\nabla\pi^*(w_\delta)
- \nabla \pi^*(w_0)\big)$.
Thus,
\begin{equation*}
\begin{aligned}
e^{\frac {rt}{\delta}}\big(z_{\delta}(t) - z_0\big)=
\frac r \delta \int_0^t e^{\frac {rs}{\delta}} \big(\nabla\pi^*(w_\delta(s))- \nabla \pi^*(w_0)\big)
{\rm d s}
\end{aligned}
\end{equation*}
and moreover,
\begin{equation}
\begin{aligned}
\label{ine:B}
\Vert z_{\delta}(t) - z_0 \Vert \le 
\frac r \delta \int_0^t L_{\pi^*} A_{\delta}(t) \frac{s^2}{2} {\rm ds} 
\le \frac {L_{\pi^*} r t^2}{6} A_{\delta}(t).
\end{aligned}
\end{equation}
If $t > \delta$, 
\begin{equation*}
\begin{aligned}
t^r \big( 
\dot z_{\delta} +\frac r t (z_{\delta} - z_0) \big) = rt^{r-1} \big(\nabla\pi^*(w_\delta)- \nabla \pi^*(w_0)\big)
\end{aligned}
\end{equation*}
i.e.,
$\frac {\rm d}{\rm dt} t^r \big( z_{\delta}(t) - z_0 \big)\!=\! rt^{r-1} \big(\nabla\pi^*(w_\delta(t)) -
\nabla \pi^*(w_0)\big).$
Similar to (\ref{ine:B}),
$\Vert z_{\delta}(t) - z_0 \Vert 
\le L_{\pi^*} r t^2 A_{\delta}(t)/6$.
Consequently,
\begin{equation}
\begin{aligned}
\label{ine:B2}
B_{\delta}(t) \le \beta_4 t A_{\delta}(t), 
~\forall t \ge 0
\end{aligned}
\end{equation}
where $\beta_4 = L_{\pi^*}r / 6$.

According to (\ref{alg:apsim}), we also obtain
\begin{equation*}
\begin{aligned}
\Vert & \dot z_{\delta}(t)\Vert 
= \frac {r}{\max(\delta, t)} \Vert \nabla \pi^*(w_{\delta}(t)) - z_{\delta}(t)\Vert \\
&\le  \frac {r}{\max(\delta, t)} \big(
\Vert \nabla \pi^*(w_{\delta}(t))-\nabla \pi^*(w_0) \Vert + \Vert  z_{\delta}(t)
- z_0\Vert \big) \\
& \le \frac{r}{\max(\delta, t)}\Big(\frac {L_{\pi^*}t^2}2 A_{\delta}(t)+ tB_{\delta}(t) \Big).
\end{aligned}
\end{equation*}
Therefore,  $C_\delta(t) \le \beta_5 t A_{\delta}(t)$, where $\beta_5 = r(L_{\pi^*}/2 + \beta_4)$.

{\color{blue}We next show $A_{\delta}, B_{\delta}$ and $C_{\delta}$ are bounded on $[0, t_m]$ for some $t_m > 0$.}
In light of (\ref{A_ineq}),
\begin{equation*}
\begin{aligned}
r A_{\delta}(t) \le \beta_0 
+ \beta_1\beta_4 t^2 A_{\delta}(t)
+ \beta_2\beta_5 t^2 A_{\delta}(t)
+ \beta_3 t^2 A_\delta(t).
\end{aligned}
\end{equation*}
Take $t_1$ such that 
$\xi_0 = r - (\beta_3 + \beta_2 \beta_5 + \beta_1 \beta_4) t_1^2 > 0$, and $t_2 = \min \{t_0, t_1\}$. 
{\color{blue}Then $A_\delta(t_2)\le \beta_0 / \xi_0$. Additionally, 
$B_{\delta}(t_2) \le \beta_0 \beta_4 t_2 / \xi_0$, and $C_{\delta}(t_2)\le \beta_0 \beta_5 t_2/ \xi_0$.

Thus, 
$\Vert \dot w_{\delta}(t)\Vert \le  \beta_0 t_2/ \xi_0$,
and $\Vert \dot z_{\delta}(t)\Vert \le  \beta_0 \beta_5 t_2/ \xi_0$ for $t \in [0, t_2]$.
In conclusion, both $z_{\delta}(t)$ and $w_{\delta}(t)$ are equi-continuous on $[0, t_2]$. It also follows that they are uniformly bounded on the same interval.}
This completes the proof.
$\hfill\square$

\emph{Remaining Proof for Theorem 1:}
Here we show that there is a unique solution to (\ref{alg:asim}).

Suppose $(z_1, w_1)$ and $(z_2, w_2)$ be two solutions of (\ref{alg:asim}).
Let $\triangle_z = z_1 - z_2$ and  $\triangle_w = w_1 - w_2$. 
Then 
\begin{equation}
\begin{aligned}
\label{pf:uniq}
\dot \triangle_z =& \frac r t \big[ 
\nabla \pi^*(w_1) -\nabla \pi^*(w_2)-\triangle_z \big] \\
\dot \triangle_w =& \frac t r \big[\! - \nabla R(z_1) + \nabla R(z_2) - \tilde H \big(\triangle_z + \frac tr \dot \triangle_z \big) \\
&+ \nabla \pi(z_1 + \frac tr \dot z_1) - \nabla \pi(z_2 + \frac tr \dot z_2) - \triangle_w \big].
\end{aligned}
\end{equation}

Define $A(t) = \sup_{[0, t]} \Vert \dot \triangle_w(s)\Vert /s$ and $B(t) = \sup_{[0, t]} \Vert\triangle_z(s)\Vert$.
By a similar procedure as the proof of (\ref{ine:B2}), we derive
$B(t) \le L_{\pi^*} r t^2 A(t)/6$.

Due to (\ref{pf:uniq}), 
$\Vert \dot \triangle_z\Vert \le  \frac rt\big(L_{\pi^*} \Vert \triangle_w\Vert + \Vert \triangle_z\Vert \big)$.
Besides,
\begin{equation*}
\begin{aligned}
\frac r t\Vert \dot \triangle_w \Vert \le &  L_R \Vert \triangle_z \Vert + \lambda_H \Vert \triangle_z\Vert
+  \frac  {\lambda_H t} {r}\Vert \dot \triangle_z\Vert + 
L_{\pi} \Vert \triangle_z \Vert \\
&+\frac { L_{\pi} t} {r} \Vert \dot \triangle_z \Vert  +
\Vert \triangle_w \Vert.
\end{aligned}
\end{equation*}

Note that $\Vert \triangle_w \Vert \le A(t) t^2/2$. As a consequence,
$$rA(t) \le \gamma_1 B(t) + \gamma_2 t^2 A(t),$$
where $\gamma_1 = L_R + 2 \lambda_H + 2L_{\pi}$ and 
$\gamma_2 = (L_{\pi^*} L_{\pi} + L_{\pi^*} \lambda_H + 1) / 2$.
Furthermore, there exists $t_4$ such that $\xi_1 := r - \gamma_2 t_4^2 > 0$, and $A(t) \le \gamma_1 B(t) / \xi_1 \le (\gamma_1 L_{\pi^*} r t^2 A(t))/ (6\xi_1), \forall t \in [0, t_4)$.
It follows that $A(t) = B(t) = 0, ~\forall t \in [0, t_5)$ for some $t_5 > 0$.
Thus, (\ref{alg:asim}) has a unique solution on $t \in [0, t_5)$. 
It is straightforward that (\ref{alg:asim}) has a unique solution
for $t \ge t_5$, and this completes the proof.
$\hfill\square$

\bibliographystyle{IEEEtran}
\bibliography{references.bib}

\end{document}